\documentclass[11pt,a4paper]{article}

\usepackage[english]{babel}
\usepackage{amssymb,amsmath,enumerate,a4wide,amsthm}
\usepackage{colortab}
\usepackage{color}
\usepackage{epsfig}
\usepackage[hidelinks]{hyperref}

\newcommand{\bb}{\mathbb}
\newcommand{\R}{\bb R}

\newcommand{\Z}{\bb Z}

\newcommand{\ceil}[1]{\left\lceil#1\right\rceil}

\newcommand{\sm}{\setminus}

\renewcommand{\epsilon}{\varepsilon}

\newtheorem{prop}{Proposition}

\newtheorem{lemma}[prop]{Lemma}

\newtheorem{thm}[prop]{Theorem}

\newtheorem*{claim-no-number}{\sc Claim}

\newtheorem{problem}{Problem}

\begin{document}

\title{Hardness of integer linear feasibility with cyclic support}

\author{Marco Di Summa\thanks{University of Padua, Department of Mathematics ``Tullio Levi-Civita'', Italy.}}

\date{October 2026}

\maketitle

\begin{abstract}
We prove that integer linear feasibility is strongly $NP$-complete even when the constraint matrix has exactly two nonzero entries in every row and column, with its support forming a single cycle. This holds already for monotone systems (one positive and one negative coefficient in each inequality) with nonzero coefficients of absolute value $2$ or $3$ and right-hand sides in $\{-1,0,1,2\}$. No variable bounds are supplied separately from the system. For systems in which each inequality involves at most two variables, this strengthens known hardness results by restricting all numerical data to fixed finite sets. For systems in which each variable occurs in at most two inequalities, this establishes, to the best of our knowledge, the first $NP$-hardness result for feasibility. We also give negative answers, unless $P=NP$, to questions on polynomial-time solvability of integer programs with bounded coefficients and bounded totally $\Delta$-modular treewidth, and on fixed-parameter tractability of integer bound propagation with respect to coefficient and right-hand-side magnitudes.
\end{abstract}

\section{Introduction}\label{sec:introduction}

Given a matrix $A\in\mathbb Z^{m\times n}$ and a vector
$b\in\mathbb Z^m$, the integer linear feasibility problem asks whether
\begin{equation}\label{eq:general-IP}
    Ax\leq b,\qquad x\in\mathbb Z^n
\end{equation}
has a solution. It is well known that this problem is $NP$-complete; see, e.g., \cite{schrijver1986}. In this paper, we are interested in sparse
instances of \eqref{eq:general-IP}. In particular, we consider the
case in which every row of $A$ has at most two nonzero entries and
the case in which every column of $A$ has at most two nonzero
entries.

In formulation \eqref{eq:general-IP}, no lower or upper bounds are supplied separately from the system.
If a bound is present, it must be described by a row of $Ax\leq b$, and therefore counts as an occurrence of the corresponding variable. This convention is immaterial when row sparsity is considered, since a bound is a one-variable inequality, but it is essential in the case of column sparsity.

The threshold two is particularly natural for sparse linear and integer
programming. In the continuous setting, linear programs with at most
two nonzero entries per row or per column admit strongly polynomial-time
algorithms \cite{DadushKohNaturaOlverVegh2024}.
The corresponding integer feasibility problems, however, behave
very differently.

\subsection{Two nonzero entries per row}

Integer systems in which every inequality contains at most two
variables are usually referred to as integer
two-variable-per-inequality, or integer TVPI, systems.
A TVPI system is called monotone if no inequality has two nonzero
coefficients with the same sign.

Lagarias \cite{Lagarias1985} proved that integer TVPI feasibility is $NP$-complete.
Strong $NP$-completeness, even for monotone systems, follows from
the construction in the proof of
\cite[Theorem~5]{EibenGanianKnopOrdyniak2018}.
The coefficients and right-hand sides in that construction are
polynomially bounded in the size of the source instance, but are
not bounded by an absolute constant.
Moreover, one variable occurs in every inequality, so the
construction does not establish hardness under the restriction
of at most two nonzero entries per column.

A substantial literature has investigated positive results for integer TVPI systems.
When finite lower and upper bounds on the variables are available,
pseudopolynomial algorithms for feasibility are known
\cite{HochbaumMegiddoNaorTamir1993,HochbaumNaor1994,BarYehudaRawitz2001}.
Efficient algorithms are available under additional restrictions on
the coefficients. For instance, one can solve in polynomial time
TVPI systems with coefficients in $\{-1,0,1\}$
\cite{LahiriMusuvathi2005,SubramaniWojciechowski2023}
and monotone TVPI systems in which the nonzero coefficients are
powers of a fixed integer in absolute value \cite{deVlas2024}.

Tamori and Kimura \cite{TamoriKimura2019} give a
polynomial-time algorithm for integer TVPI systems when all strongly connected components
of the associated directed graph have bounded size.
This graph has one node per variable and an arc from $x_i$ to
$x_j$ for each inequality $cx_i-dx_j\le\beta$ with $c,d>0$.

\subsection{Two nonzero entries per column}
The restriction that every column of $A$ contains at most two
nonzero entries means that each variable may
occur in at most two inequalities. Again, the continuous problem is
well behaved and closely related to generalized network
optimization \cite{Hochbaum2004,DadushKohNaturaOlverVegh2024}.

Hardness results are known for closely related integer programming
models, but they do not directly apply to \eqref{eq:general-IP}.
In particular, feasibility of
\[    Ax=b,\qquad x\geq 0, \qquad x\in\Z^n\]
is strongly $NP$-complete even when $A$ has at most two nonzero entries
per column; see \cite{Hochbaum2004,PritchardChakrabarty2011}.
However, replacing each equality by two inequalities and expressing
nonnegativity by additional inequalities may produce five
occurrences per variable.

Pritchard and Chakrabarty \cite{PritchardChakrabarty2011} also discuss integer linear {\em optimization} problems in which each variable occurs in at most two inequalities, including its bounds. They note that $NP$-hardness was known in this setting, for instance via Unbounded Knapsack, while strong $NP$-hardness remained open. However, the usual formulation of Unbounded Knapsack is always feasible, and imposing a lower bound on the objective value introduces a third occurrence of each variable. Thus, this does not establish $NP$-hardness of feasibility for \eqref{eq:general-IP}.

A recent result
\cite{BrianskiLassotaPekarkovaPilipczukReuter2026}
establishes $NP$-hardness of feasibility for integer programs with a ``path-like'' constraint matrix whose entries are bounded in absolute value by
$8$. The formulation is
\[    Ax=b,\qquad \ell\leq x\leq u,\qquad x\in\Z^n\]
and the explicit variable bounds $\ell$ and $u$ play an essential role in the
model. Hence that result does not settle the problem considered here.

To the best of our knowledge, the complexity of
\eqref{eq:general-IP} under the restriction that every column of
$A$ has at most two nonzero entries has remained open, even with
respect to ordinary $NP$-hardness.

Positive results are known when sparsity is combined with
additional arithmetic structure. In particular, Fiorini et al.~\cite{FioriniJoretWeltgeYuditsky2025}
give a strongly polynomial-time algorithm to solve integer programs with at
most two nonzero entries per row where all subdeterminants are bounded in absolute value by a fixed constant. They
also observe that the analogous case with at most two nonzero
entries per column can be solved in strongly polynomial time, via a
reduction to $b$-matching. Thus, bounded subdeterminants yield
tractability for both sparsity patterns, whereas bounding the
individual entries of the matrix does not, as our result below shows.

\subsection{Our result}

We prove that integer feasibility is $NP$-complete even for square
monotone systems whose associated digraph consists of a single
cycle and whose coefficients and right-hand sides belong to fixed
finite sets. In particular, the matrices used in our reduction
have exactly two nonzero entries in every row and column.

\begin{problem}[{\normalfont\scshape Cyclic-IP-Feasibility}]\label{problem1}
\leavevmode\\
{\normalfont\scshape Input:}
An integer $n\geq 2$ and integers $p_i,q_i,b_i$ for
$i=1,\dots,n$, with $p_i,q_i>0$ for every $i$.\\
{\normalfont\scshape Question:}
Decide whether there exists $x\in\mathbb Z^n$ satisfying
\begin{equation}\label{eq:system}
    p_i x_i-q_i x_{i+1}\leq b_i,
    \qquad i=1,\dots,n,
\end{equation}
where $x_{n+1}=x_1$.
\end{problem}

Our main result, which we prove in Section \ref{sec:proof}, is the following.

\begin{thm}\label{thm:main}
Problem {\normalfont\scshape Cyclic-IP-Feasibility} is $NP$-complete even under the additional assumption
that 
\[ (p_i,q_i)\in\{(2,3),(3,2)\}   \qquad\text{and}\qquad
    b_i\in\{-1,0,1,2\}
\qquad \mbox{for }i=1,\dots,n.\]
Therefore, Problem {\normalfont\scshape Cyclic-IP-Feasibility} is strongly $NP$-complete.
\end{thm}

Viewed as a result for integer TVPI systems,
Theorem \ref{thm:main} strengthens the strong $NP$-completeness
result of Eiben et al.~\cite{EibenGanianKnopOrdyniak2018}
by restricting the nonzero coefficients to absolute values
$2$ and $3$ and the right-hand sides to $\{-1,0,1,2\}$.
For any fixed integer $a\ge2$, feasibility of monotone TVPI systems
with nonzero coefficients of absolute value in $\{1,a\}$ can be
checked in polynomial time \cite{deVlas2024}.
Thus, unless $P=NP$, two distinct absolute values larger than $1$
are necessary for $NP$-hardness of monotone TVPI feasibility.
Our choice of $2$ and $3$ is therefore optimal both in the number
of distinct absolute values and in their magnitude.

Viewed by columns, Theorem \ref{thm:main} establishes strong
$NP$-completeness of feasibility for integer programs in which
each variable occurs in at most two inequalities.
Thus, it resolves the question discussed above, for which even
ordinary $NP$-hardness of feasibility was not known.
It also settles the question of strong $NP$-hardness for
optimization raised in \cite{PritchardChakrabarty2011}. 

The cyclic structure also gives a sharp contrast with
\cite[Theorem~7]{EibenGanianKnopOrdyniak2018}, which shows that
integer feasibility is polynomial-time solvable when the variable-constraint incidence
graph is acyclic. This is the bipartite graph
with one vertex for each variable and each inequality, and an edge
for each nonzero entry of $A$.
For \eqref{eq:system}, this graph is a single cycle.
Thus, Theorem \ref{thm:main} shows that allowing one cycle
already yields strong $NP$-completeness, even when all coefficients
and right-hand sides are bounded by a constant.

We also note that, for systems in integer variables whose variable-constraint incidence graph is a single cycle, the monotone case captures all potentially infeasible instances, up to sign changes of the variables. Indeed, if some column contains two nonzero entries with the same sign, the system is always feasible: the inequalities not involving the corresponding variable form a path and can be satisfied successively, after which that variable can be chosen with sufficiently large absolute value and the appropriate sign to satisfy its two incident inequalities. Otherwise, suitable sign changes of the variables make the system monotone.

Finally, Theorem \ref{thm:main} gives negative answers (unless $P=NP$) to questions on
poly\-nomial-time solvability of integer programs with bounded coefficients and bounded totally $\Delta$-modular treewidth \cite{McFarland2026}, and on fixed-parameter
tractability of integer bound propagation with respect to coefficient and right-hand-side magnitudes
\cite{BordeauxKatsirelosNarodytskaVardi2011}. These consequences are discussed in Section~\ref{sec:consequences}.

\subsection{Proof idea}

We prove Theorem~\ref{thm:main} by a reduction from the following
not-all-equal satisfiability problem, which is known to be $NP$-complete \cite{Schaefer1978}.

\begin{problem}[{\normalfont\scshape Monotone NAE-3SAT}]
\leavevmode\\
{\normalfont\scshape Input:}
A nonempty collection of clauses $S_1,\dots,S_m$, each consisting of three distinct unnegated Boolean variables.\\
{\normalfont\scshape Question:}
Decide whether the variables occurring in the clauses can be assigned Boolean values such that no clause contains three variables with the same value.
\end{problem}

The term ``monotone'' indicates that no clause contains a negated variable.
We denote the variables occurring in the clauses by
$v_0,\dots,v_{N-1}$ and identify each clause $S_k$ with the corresponding three-element subset of $\{0,\dots,N-1\}$.
 %We will often identify the Boolean values true and false with the bits 1 and 0, respectively.

We fix an instance of \textsc{Monotone NAE-3SAT} and associate with each integer $z$ a Boolean assignment, in such a way that every assignment is represented by some integer. We then construct an instance of \textsc{Cyclic-IP-Feasibility}, which is conveniently described in terms of \emph{transfer maps}.
For integer variables $x_i,x_{i+1}$, an inequality
$p_i x_i-q_i x_{i+1}\le b_i$, with $p_i,q_i>0$,
is equivalent to $x_{i+1}\ge F_i(x_i)$, where
\[
F_i(z):=\ceil{\frac{p_i z-b_i}{q_i}}
\]
is called the transfer map of the inequality.
For a \emph{block} of inequalities $p_i x_i-q_i x_{i+1}\le b_i$,
$i=1,\dots,\ell-1$, the transfer map is
$G:=F_{\ell-1}\circ\cdots\circ F_1$. Since each $F_i$ is nondecreasing,
$G(z)$ is the minimum possible value of $x_\ell$ when $x_1=z$. 
These transfer maps are one-dimensional affine Chv\'atal
functions, and their construction amounts to eliminating the
intermediate integer variables; see
\cite{BlairJeroslow1982,BasuMartinRyanWang2019}.

The key step is to construct, for each clause $S$, a block whose
transfer map $C_S$ has the property that $C_S(z)\ge z$ for every integer $z$,
with equality if and only if the assignment encoded by $z$
satisfies $S$. Each block contains $O(N)$ inequalities and uses only
the coefficients and right-hand sides allowed in
Theorem \ref{thm:main}.

Finally, we concatenate the clause blocks and identify the last
variable with the first, obtaining a cyclic system. As we prove below, a cyclic system of this form is feasible
if and only if its transfer map $F$ satisfies $F(z)\le z$ for some
integer $z$. Since the transfer map of each clause block satisfies $C_S(z)\ge z$, this
condition forces all these transfer maps to fix this particular $z$. By construction, this is equivalent to every clause
being satisfied by the assignment encoded by $z$. As every Boolean
assignment has such an encoding, the reduction follows.

\section{Detailed proof}\label{sec:proof}

We first prove the cyclic feasibility criterion stated in the
introduction.

\begin{lemma}\label{lem:cyclic-transfer}
For $i=1,\dots,n$, let $p_i,q_i,b_i$ be as in Problem \ref{problem1}. Let $F_i$ be the transfer map of the $i$th inequality in \eqref{eq:system}, i.e.,
$F_i(z):=\ceil{(p_i z-b_i)/q_i}$, and set 
$F:=F_n\circ\cdots\circ F_1$.
System \eqref{eq:system} has an integer solution if and only if
there exists $z\in\mathbb Z$ such that $F(z)\le z$.
\end{lemma}

\begin{proof}
Since $x_2\ge F_1(x_1)$ in every feasible solution and $x_2$
has a negative coefficient only in the first inequality of
\eqref{eq:system}, replacing $x_2$ by $F_1(x_1)$ preserves
feasibility. Repeating this operation for $x_3,\dots,x_n$ shows that
system \eqref{eq:system} has an integer solution if and only if there exists $x_1\in\Z$ such that, defining
$x_{i+1}=F_i(x_i)$ for $i=1,\dots,n-1$, $x_1$ and $x_n$
satisfy the last inequality of \eqref{eq:system}.
The last condition is equivalent to
$x_1\ge F_n(x_n)=F(x_1)$.
%
%Suppose that $x\in\mathbb Z^n$ satisfies \eqref{eq:system}.
%Then $x_{i+1}\ge F_i(x_i)$ for $i=1,\dots,n$, where
%$x_{n+1}=x_1$. Since each $F_i$ is nondecreasing,
%\[x_1\ge F_n(x_n)\ge F_n(F_{n-1}(x_{n-1}))
%\ge\cdots\ge F_n(F_{n-1}(\cdots F_1(x_1)))= F(x_1).\]
%Thus $z:=x_1$ satisfies $F(z)\le z$.
%
%Conversely, suppose that $F(z)\le z$ for some $z\in\mathbb Z$.
%Set $x_1:=z$ and, recursively, $x_{i+1}:=F_i(x_i)$ for
%$i=1,\dots,n-1$. The first $n-1$ inequalities in \eqref{eq:system} hold by
%construction, and the last one follows from
%$F_n(x_n)=F(x_1)\le x_1$.
\end{proof}

For $a,c\in\mathbb Z$, define
\[
U_a(z):=\ceil{\frac{3z-a}{2}},
\qquad
V_c(z):=\ceil{\frac{2z-c}{3}}.
\]
These are transfer maps of  inequalities of the form $3x_i-2x_{i+1}\le a$ and
$2x_i-3x_{i+1}\le c$, respectively.

A sequence $T_1,\dots,T_\ell$ of maps of this type specifies a block
with integer variables $y_0,\dots,y_\ell$ and inequalities
$y_r\ge T_r(y_{r-1})$, $r=1,\dots,\ell$, each written in the
corresponding linear form. The \emph{input} and \emph{output} variables
are $y_0$ and $y_\ell$, respectively, and the transfer map of the block is $T_\ell\circ\cdots\circ T_1$.

Fix an instance of \textsc{Monotone NAE-3SAT} with variables $v_0,\dots,v_{N-1}$ and clauses $S_1,\dots,S_m$.
We identify the Boolean values false and true with the bits $0$ and $1$, respectively.

%\subsection{Encoding Boolean assignments}

Set $U:=U_0$. For $z\in\mathbb Z$ and an integer $j\ge0$,
let
\begin{equation}\label{eq:zj}
z_j:=U^j(z),\qquad \epsilon_j(z):=z_j\bmod 2\in\{0,1\},
\end{equation}
where $U^j$ denotes the $j$-fold composition
of $U$, with $U^0$ being the identity.
The integer $z$ encodes the assignment that gives $v_j$ the value $\epsilon_j(z)$ for $j=0,\dots,N-1$. The next lemma shows
that every Boolean assignment is represented in this way.
We denote by $[z]_{2^H}$ the residue class of $z$ modulo $2^H$.

\begin{lemma}\label{lem:parity-encoding}
For every integer $H\geq 1$, the map
\[
\Psi_H:\mathbb Z/2^H\mathbb Z\longrightarrow\{0,1\}^H,
\qquad
[z]_{2^H}\longmapsto
\bigl(\epsilon_0(z),\dots,\epsilon_{H-1}(z)\bigr)
\]
is a bijection.
\end{lemma}

\begin{proof}
We first verify that $\Psi_H$ is well defined. For every integer
$k\geq 1$ and $t\in\mathbb Z$,
\[U(z+2^k t)=\left\lceil\frac{3z}{2}+3\cdot 2^{k-1}t\right\rceil=U(z)+3\cdot 2^{k-1}t.\]
Consequently,
\[z\equiv z'\pmod{2^k}
\quad\Longrightarrow\quad
U(z)\equiv U(z')\pmod{2^{k-1}}.\]
By iterating this implication, we obtain 
\[z\equiv z'\pmod{2^H} \quad\Longrightarrow\quad
U^j(z)\equiv U^j(z')\pmod{2^{H-j}}
\quad \text{for }j=0,\dots,H-1.\]
In particular, recalling \eqref{eq:zj}, $\epsilon_j(z)=\epsilon_j(z')$ for
$j=0,\dots,H-1$. Hence $\Psi_H$ depends only on the residue class
of $z$ modulo $2^H$.

Since both $\mathbb Z/2^H\mathbb Z$ and $\{0,1\}^H$ have
cardinality $2^H$, to prove that $\Psi_H$ is bijective it is sufficient to show injectivity. Suppose that $z,z'\in\mathbb Z$ satisfy
\begin{equation}\label{eq:eps}
\epsilon_j(z)=\epsilon_j(z')
\qquad
\text{for }j=0,\dots,H-1.
\end{equation}
Similar to \eqref{eq:zj}, set $z'_j:=U^j(z')$ for every $j\ge0$.
Since $z_{j+1}=U(z_j)=\lceil{3z_j/2}\rceil$ and $z'_{j+1}=U(z'_j)=\lceil{3z'_j/2}\rceil$, we have
\begin{equation}\label{eq:z-eps}
2z_{j+1}=3z_j+\epsilon_j(z),
\qquad
2z'_{j+1}=3z'_j+\epsilon_j(z').
\end{equation}
Using \eqref{eq:eps}, we obtain
$2(z_{j+1}-z'_{j+1})=3(z_j-z'_j)$ for $j=0,\dots,H-1$, which implies $2^H(z_H-z'_H)=3^H(z-z')$.
Therefore $2^H$ divides $3^H(z-z')$. Since $2^H$ and $3^H$
are coprime, it follows that $2^H$ divides $z-z'$.
Thus $z\equiv z'\pmod{2^H}$, and $\Psi_H$ is injective.
\end{proof}

We now fix a clause 
\[
S\subseteq\{0,\dots,N-1\},
\qquad |S|=3,
\]
and set $H:=1+\max S$.
In particular, $H-1\in S$.

A not-all-equal clause is violated in exactly two cases: all the bits indexed
by $S$ are equal to $1$, or all of them are equal to $0$.
Accordingly, we construct two \emph{penalty blocks}. The transfer map of the first block adds one unit to its input precisely in the all-one case, while that of the second block adds one unit
precisely in the all-zero case.

%\medskip
%\noindent
%\textbf{The all-one penalty block.}

We start with the all-one penalty block.
For $j=0,\dots,H-1$, define
\[
c_j:=
\begin{cases}
0 & \mbox{if }j=H-1,\\
2 & \mbox{if }j\in S\setminus\{H-1\},\\
1 & \mbox{if }j\notin S.
\end{cases}
\]
Let $\mathcal B_S^1$ be the block of inequalities corresponding to the sequence of maps
\[
\underbrace{U,\dots,U}_{H\text{ times}},
V_{c_{H-1}},V_{c_{H-2}},\dots,V_{c_0}.
\]
Its transfer map is
\begin{equation}\label{eq:BS1}
B_S^1:=V_{c_0}\circ V_{c_1}\circ\cdots\circ V_{c_{H-1}}\circ U^H.
\end{equation}

\begin{lemma}\label{lem:all-one-block}
For every $z\in\mathbb Z$,
\[B_S^1(z)=z+\min\{\epsilon_j(z):j\in S\}.\]
\end{lemma}

\begin{proof}
Using the definition of $V_c$ and the first equality in \eqref{eq:z-eps},
for every $\delta\in\{0,1\}$, $c\in\{0,1,2\}$ and $j\ge0$ we obtain
\begin{equation}\label{eq:down-up-identity}
V_c(z_{j+1}+\delta)=\ceil{\frac{2z_{j+1}+2\delta-c}{3}}
=z_j+\ceil{\frac{\epsilon_j+2\delta-c}{3}},
\end{equation}
where we write $\epsilon_j$ instead of $\epsilon_j(z)$ to simplify the notation. A direct verification of the four possibilities
$(\epsilon_j,\delta)\in\{0,1\}^2$ shows that
\begin{equation}\label{eq:boolean-rules-one}
\ceil{\frac{\epsilon_j+2\delta-c}{3}}
=
\begin{cases}
\max\{\epsilon_j,\delta\} & \mbox{if }c=0,\\
\delta                & \mbox{if }c=1,\\
\min\{\epsilon_j,\delta\} & \mbox{if }c=2.
\end{cases}
\end{equation}

Starting from $w_H:=z_H=U^H(z)$, define
\[
w_j:=V_{c_j}(w_{j+1}),
\qquad
j=H-1,H-2,\dots,0.
\]
Then, by \eqref{eq:BS1}, $B_S^1(z)=w_0$.

Define $\delta_j=w_j-z_j$ for every $j\in\{0,\dots,H\}$.
Since $w_H=z_H$, $\delta_H=0$.
For $j=H-1,\dots,0$, assuming that $\delta_{j+1}\in\{0,1\}$,
we obtain
\[
w_j=V_{c_j}(w_{j+1})
=V_{c_j}(z_{j+1}+\delta_{j+1})
=z_j+\ceil{\frac{\epsilon_j+2\delta_{j+1}-c_j}{3}},
\]
where the last equality follows from \eqref{eq:down-up-identity}.
Therefore, using \eqref{eq:boolean-rules-one},
\[
\delta_j=
\begin{cases}
\max\{\epsilon_j,\delta_{j+1}\} & \mbox{if }c_j=0,\\
\delta_{j+1} & \mbox{if }c_j=1,\\
\min\{\epsilon_j,\delta_{j+1}\} & \mbox{if }c_j=2.
\end{cases}
\]
It follows by backward induction that $\delta_j\in\{0,1\}$
for $j=H,H-1,\dots,0$.
Furthermore, since $c_{H-1}=0$ and $\delta_H=0$,
we have $\delta_{H-1}=\epsilon_{H-1}$.
For $0\le j\le H-2$, $c_j\in\{1,2\}$, and therefore
$\delta_0=1$ if and only if $\delta_{H-1}=1$ and
$\epsilon_j=1$ whenever $c_j=2$, i.e., whenever
$j\in S\sm\{H-1\}$.
Since $\delta_{H-1}=\epsilon_{H-1}$ and $H-1\in S$,
it follows that $\delta_0=1$ if and only if $\epsilon_j=1$
for every $j\in S$.
As $B_S^1(z)=w_0=z+\delta_0$, this concludes the proof.
\end{proof}

% \medskip
% \noindent
% \textbf{The all-zero penalty block.}

We next construct the all-zero penalty block. For
$j=0,\dots,H-1$, define
\[
a_j:=
\begin{cases}
-1 &\mbox{if } j=H-1,\\
1 &\mbox{if } j\in S\setminus\{H-1\},\\
0 &\mbox{if } j\notin S.
\end{cases}
\]
Let $\mathcal B_S^0$ be the block of inequalities corresponding to the sequence of maps
\[
U_{a_0},U_{a_1},\dots,U_{a_{H-1}},
\underbrace{V_1,\dots,V_1}_{H\text{ times}}.
\]
Its transfer map is
$B_S^0:=V_1^H\circ U_{a_{H-1}}\circ\cdots\circ U_{a_1}\circ U_{a_0}$.

\begin{lemma}\label{lem:all-zero-block}
For every $z\in\mathbb Z$,
\[B_S^0(z)=z+\min\{1-\epsilon_j(z):j\in S\}.\]
\end{lemma}

\begin{proof}
We first establish some identities. For $z\in\mathbb Z$ and $\epsilon:=\epsilon_0(z)$, 
it is straightforward to verify that
\[U_{-1}(z)=\frac{3z+2-\epsilon}{2},\qquad
U_0(z)=\frac{3z+\epsilon}{2},\qquad
U_1(z)=\frac{3z-\epsilon}{2}.\]
For $\delta\in\{0,1\}$, the analysis of the four possibilities
$(\epsilon,\delta)\in\{0,1\}^2$ then gives
\begin{align}
V_1(U_{-1}(z)+\delta)
&=\ceil{\frac{3z+1-\epsilon+2\delta}{3}}
=z+\max\{1-\epsilon,\delta\},
\label{eq:zero-init}\\
V_1(U_0(z)+\delta)
&=\ceil{\frac{3z+\epsilon+2\delta-1}{3}}
=z+\delta,
\label{eq:zero-pass}\\
V_1(U_1(z)+\delta)
&=\ceil{\frac{3z-\epsilon+2\delta-1}{3}}
=z+\min\{1-\epsilon,\delta\}.
\label{eq:zero-and}
\end{align}

Define
\begin{equation}\label{eq:R-init}
R_{H-1}:=V_1\circ U_{a_{H-1}}=V_1\circ U_{-1},
\end{equation}
and, recursively for $i=H-2,\dots,0$,
\begin{equation}\label{eq:Ri}
R_i:=V_1\circ R_{i+1}\circ U_{a_i}.
\end{equation}
Expanding the recursion gives
$R_0=V_1^H\circ U_{a_{H-1}}\circ\cdots\circ U_{a_0}=B_S^0$.

We claim that, for every $i=0,\dots,H-1$ and every
$z\in\mathbb Z$,
\begin{equation}\label{eq:Rj-invariant}
R_i(z)=z+\min \{1-\epsilon_{j-i}(z):j\in S,\,j\ge i\}.
\end{equation}
As $R_0=B_S^0$, for $i=0$ this claim proves the lemma.

We prove the claim by backward induction on $i$. For $i=H-1$,
recall that $H-1=\max S$. By \eqref{eq:R-init} and \eqref{eq:zero-init} with
$\delta=0$,
\[R_{H-1}(z)=V_1(U_{-1}(z))=z+1-\epsilon_0(z).\]
This is precisely \eqref{eq:Rj-invariant} for $i=H-1$.

Suppose now that \eqref{eq:Rj-invariant} holds for the index $i+1$,
where $0\le i\le H-2$. Applying the induction hypothesis with
$U_{a_i}(z)$ in place of $z$ gives
\begin{equation}\label{eq:all0-R}
R_{i+1}(U_{a_i}(z))=U_{a_i}(z)+\delta,
\end{equation}
where
\begin{equation}\label{eq:all0-delta}
\delta=\min\left\{1-\epsilon_{j-i-1}(U_{a_i}(z)):
j\in S,\,j\ge i+1\right\}\in\{0,1\}.
\end{equation}

We now distinguish two cases.

\begin{itemize}
\item 
Assume that $i\notin S$. In this case $a_i=0$. Using \eqref{eq:Ri}, \eqref{eq:all0-R}, and then \eqref{eq:zero-pass}, we obtain
\begin{equation}\label{eq:Ri-chain}
R_i(z)=V_1(R_{i+1}(U_0(z)))=V_1(U_0(z)+\delta)=z+\delta.
\end{equation}
Since $U_{a_i}(z)=U_0(z)=U(z)$ and $\epsilon_{j-i}(z)=\epsilon_{j-i-1}(U(z))$, \eqref{eq:all0-delta} gives
\begin{equation}\label{eq:delta-chain}
\begin{split}
\delta&=\min\{1-\epsilon_{j-i}(z):j\in S,\,j\ge i+1\}\\
&=\min\{1-\epsilon_{j-i}(z):j\in S,\,j\ge i\},
\end{split}
\end{equation}
where the last equality follows from the assumption $i\notin S$.
Together, \eqref{eq:Ri-chain} and \eqref{eq:delta-chain} imply \eqref{eq:Rj-invariant}.
\item 
Suppose now that $i\in S$. In this case $a_i=1$. Using \eqref{eq:Ri}, \eqref{eq:all0-R}, and then \eqref{eq:zero-and}, we obtain
\begin{equation}\label{eq:Ri-case2}
R_i(z)=V_1(R_{i+1}(U_1(z)))=V_1(U_1(z)+\delta)=z+\min\{1-\epsilon_0(z),\delta\}.
\end{equation}

If $\epsilon_0(z)=1$, then
$R_i(z)=z$, as $\delta\in\{0,1\}$ by \eqref{eq:all0-delta}. Since the minimum in \eqref{eq:Rj-invariant}
contains the term $1-\epsilon_0(z)=0$, the claimed identity holds.

If $\epsilon_0(z)=0$, then $z$ is even, so
$U_{a_i}(z)=U_1(z)=U(z)$. Thus \eqref{eq:all0-delta} gives
\[\delta=\min\{1-\epsilon_{j-i}(z):j\in S,\,j\ge i+1\}.\]
Since $1-\epsilon_0(z)=1$ and $\delta\in\{0,1\}$, \eqref{eq:Ri-case2} and the above equality yield
\[
R_i(z)=z+\delta
=z+\min\{1-\epsilon_{j-i}(z):j\in S,\,j\ge i\},
\]
which is precisely \eqref{eq:Rj-invariant}.
\end{itemize}

This completes the proof of the lemma.
\end{proof}

% \medskip
% \noindent
% \textbf{The clause block.}

Let $\mathcal C_S$ be the block obtained by placing
$\mathcal B_S^1$ first and $\mathcal B_S^0$ second. Its transfer
map is
$C_S:=B_S^0\circ B_S^1$.

\begin{lemma}\label{lem:clause-block}
For every $z\in\mathbb Z$,
$C_S(z)\geq z$.
Moreover, $C_S(z)=z$
if and only if the set $\{\epsilon_j(z):j\in S\}$ contains both 0 and 1.
\end{lemma}

\begin{proof}
By Lemmas \ref{lem:all-one-block} and \ref{lem:all-zero-block}, $B_S^1(w)\geq w$ and $B_S^0(w)\geq w$
for every $w\in\mathbb Z$; thus, $C_S(z)=B_S^0(B_S^1(z))\geq z$.

If $\{\epsilon_j(z):j\in S\}$ contains both 0 and 1, then, again by Lemmas \ref{lem:all-one-block} and \ref{lem:all-zero-block}, $B_S^1(z)=B_S^0(z)=z$, and therefore $C_S(z)=z$.

If $\epsilon_j(z)=0$ for all $j\in S$, then $B_S^1(z)=z$ and $B_S^0(z)=z+1$, so $C_S(z)=z+1$.

Finally, if $\epsilon_j(z)=1$ for all $j\in S$, then $B_S^1(z)=z+1$, and since $B_S^0(w)\geq w$ for every $w\in\Z$, $C_S(z)=B_S^0(z+1)\geq z+1$.
\end{proof}

% \medskip
% \noindent
% \textbf{Construction of the cyclic instance.}

We now have all the ingredients needed to prove the main theorem.

\begin{proof}[Proof of Theorem \ref{thm:main}]
Membership in $NP$ follows from the standard fact that integer linear feasibility is in $NP$; see, e.g., \cite{Papadimitriou1981}.

To prove strong $NP$-hardness, we construct a cyclic system from
the instance of \textsc{Monotone NAE-3SAT} fixed above.
For each clause $S_k$, we construct the block
$\mathcal C_{S_k}$, and then concatenate these blocks in the order
$\mathcal C_{S_1},\mathcal C_{S_2},\dots,\mathcal C_{S_m}$.
Finally, we identify the output variable of the last block with the
input variable of the first block. This produces a cyclic system of
the form \eqref{eq:system}.

Let $F:=
C_{S_m}\circ\cdots\circ C_{S_1}$
be the transfer map of the complete sequence of clause blocks.
By Lemma \ref{lem:clause-block},
$C_{S_k}(w)\geq w$
for every $k\in\{1,\dots,m\}$ and every $w\in\mathbb Z$. Consequently,
\begin{equation}\label{eq:augmenting}
F(z)\geq z
\qquad
\text{for every }z\in\mathbb Z.
\end{equation}

We show that the formula is NAE-satisfiable if and only if the
constructed cyclic system is feasible.
Suppose first that the formula has a satisfying assignment
$\alpha:\{v_0,\dots,v_{N-1}\}\longrightarrow\{0,1\}$.
By Lemma \ref{lem:parity-encoding}, there exists $z\in\mathbb Z$
such that
$\epsilon_j(z)=\alpha(v_j)$
for $j=0,\dots,N-1$.
Since every clause contains both a variable assigned $0$ and a variable
assigned $1$, Lemma \ref{lem:clause-block} implies
$C_{S_k}(z)=z$ for $k=1,\dots,m$, and thus $F(z)=z$.
By Lemma \ref{lem:cyclic-transfer}, the constructed cyclic system
has an integer solution.

Conversely, suppose that the constructed cyclic system has an
integer solution. By Lemma \ref{lem:cyclic-transfer}, there exists
$z\in\mathbb Z$ such that $F(z)\leq z$.
Together with \eqref{eq:augmenting}, this gives $F(z)=z$.
Define
\[\zeta_0:=z\qquad\mbox{and}\qquad \zeta_k:=C_{S_k}(\zeta_{k-1}),
\quad k=1,\dots,m.\]
Lemma \ref{lem:clause-block} gives $\zeta_0\leq \zeta_1\leq\cdots\leq \zeta_m$.
On the other hand,
$\zeta_m=F(\zeta_0)=\zeta_0$, so $\zeta_0=\zeta_1=\cdots=\zeta_m$.
In particular, $C_{S_k}(z)=z$ for every $k$.
By Lemma \ref{lem:clause-block}, each set $\{\epsilon_j(z):j\in S_k\}$ contains both $0$ and $1$. Hence, the assignment
$\alpha(v_j):=\epsilon_j(z)$ for  $j=0,\dots,N-1$
satisfies every NAE clause.

% \medskip
% \noindent
% \textbf{Size and numerical values of the constructed instance.}

We finally discuss the size and numerical values of the constructed instance.
For a fixed clause $S$, each of the blocks $\mathcal B_S^1,\mathcal B_S^0$
contains $H$ maps of type $U_a$ and $H$ maps of type $V_c$,
where $H=1+\max S\le N$.
Consequently, the clause block $\mathcal C_S$ contains at most
$4N$ inequalities. All variables introduced in the construction
are distinct except for the identifications used to concatenate
the blocks and close the cycle. Thus every variable occurs in
exactly two inequalities, and the complete cyclic system contains
at most $4Nm$ variables and inequalities.
Since the inequalities are specified explicitly, the reduction can be carried out in
polynomial time.

Every inequality in the construction is of one of the forms
\[3x_i-2x_{i+1}\leq a \quad\mbox{with } a\in\{-1,0,1\}
\qquad\mbox{or}\qquad
2x_i-3x_{i+1}\leq c\quad\mbox{with } c\in\{0,1,2\}.\]
Thus all right-hand sides belong to $\{-1,0,1,2\}$.

Since the reduction is polynomial and produces only coefficients and right-hand sides whose absolute values are bounded by $3$, it proves strong $NP$-hardness. This completes the proof.
\end{proof}

We remark that each of the blocks $\mathcal B_S^1$ and
$\mathcal B_S^0$ contains the same number of maps of types $U_a$ and $V_c$. Therefore, in the complete cyclic system,
the pairs $(2,3)$ and $(3,2)$
occur equally many times. 
Then, if $A$ denotes the matrix of the constructed cyclic system, we have
\[\det A=\prod_{i=1}^n p_i-\prod_{i=1}^n q_i=0.\]
This is not surprising, as a system of inequalities with a square nonsingular constraint matrix defines a translated full-dimensional pointed cone, which always contains an integer point.
Thus, for cyclic systems of inequalities, integer feasibility can be nontrivial only if the constraint matrix is singular.

Furthermore, deleting the last row and column of $A$ yields an $(n-1)\times(n-1)$ submatrix with determinant $\prod_{i=1}^{n-1}p_i$, whose absolute value is at least $2^{n-1}$. Thus, although the entries of $A$ are bounded by a constant,
there are subdeterminants whose absolute values grow
exponentially with $n$. Therefore, there is no contradiction between Theorem \ref{thm:main} and the polynomial-time algorithms given in \cite{FioriniJoretWeltgeYuditsky2025}, where the boundedness of all subdeterminants is crucial.

\section{Further consequences}\label{sec:consequences}

We show how Theorem \ref{thm:main} implies negative answers to two additional open questions.

\subsection{A consequence for bounded TDM-treewidth}

McFarland \cite{McFarland2026} introduces a parameter called
totally $\Delta$-modular treewidth (TDM-treewidth) for matrices
with two nonzero entries per row.
Let $G(A)$ be the graph with one vertex for each column of $A$
and one edge for each row, joining the columns containing its
two nonzero entries.
When every column contains an entry of absolute value greater
than $1$, the TDM-treewidth equals the treewidth of $G(A)$ plus
one; this follows directly from
\cite[Definition~1.2]{McFarland2026}.

McFarland \cite[Question~7.2]{McFarland2026} asks whether integer
programs can be solved in polynomial time when both the absolute values of the matrix
entries and the TDM-treewidth are bounded by constants, without
requiring finite bounds on the variables.

Theorem \ref{thm:main} gives a negative answer unless $P=NP$,
already for feasibility.
Indeed, every column of a matrix produced by our reduction
contains entries of absolute value greater than $1$, and $G(A)$
is a cycle, which has treewidth two.
Thus, the matrix entries have absolute value at most $3$ and
the TDM-treewidth is exactly $3$.

\subsection{A consequence for bound propagation}

Suppose that each integer variable is initially restricted to a
bounded interval. Bound propagation repeatedly tightens these intervals using
one inequality together with the current variable bounds. For example, the
inequality $2x_i-3x_j\le1$ and the current bound $x_i\ge4$
imply $x_j\ge3$.
The process stops when an interval becomes empty or no bound
can be tightened further. The associated decision problem asks
whether all intervals are nonempty at termination.

Bordeaux et al.~\cite[Section~3.4.3]{BordeauxKatsirelosNarodytskaVardi2011}
raise the question of whether bounding the absolute values of
the coefficients and right-hand sides yields fixed-parameter
tractability. Writing $\beta$ for this bound and $L$ for the
binary input size, this would mean an algorithm with running
time $f(\beta)L^{O(1)}$.
The endpoints of the initial intervals are also given as input,
in binary, but their absolute values need not be bounded by
$\beta$. In particular, such an algorithm would run in polynomial
time for each fixed $\beta$.

Theorem \ref{thm:main} rules this out unless $P=NP$.
To see this, consider an instance satisfying the restrictions
of the theorem. Standard bounds on the size of integer solutions
\cite{Papadimitriou1981} allow us to compute in polynomial time
an integer $B$ such that restricting every variable to $[-B,B]$
preserves feasibility.
We use these intervals as the initial domains for bound
propagation, while the coefficients and right-hand sides of
the inequalities remain bounded in absolute value by $3$.

For the cyclic systems considered here, propagation terminates
with nonempty intervals if and only if the system has an integer
solution within the initial intervals.
Indeed, a feasible solution can never be removed by propagation.
Conversely, suppose that propagation terminates with nonempty
intervals, and let $\ell_i$ be their lower endpoints.
The inequality $p_i x_i-q_i x_{i+1}\le b_i$, together with
$x_i\ge\ell_i$, implies $x_{i+1}\ge F_i(\ell_i)$.
Since no lower bound can be tightened further, we must have
$\ell_{i+1}\ge F_i(\ell_i)$ for every $i$, with indices taken
cyclically. Thus, the vector of lower endpoints is feasible.
Together with Theorem \ref{thm:main} (and membership in $NP$), this shows that the decision
problem is $NP$-complete even for cyclic systems with $\beta=3$.

\bibliographystyle{plain}
\bibliography{bibliography}

@article{Papadimitriou1981,
  author  = {Papadimitriou, Christos H.},
  title   = {On the Complexity of Integer Programming},
  journal = {Journal of the ACM},
  volume  = {28},
  number  = {4},
  pages   = {765--768},
  year    = {1981},
  doi     = {10.1145/322276.322287}
}

@article{Lagarias1985,
  author  = {Lagarias, Jeffrey C.},
  title   = {The Computational Complexity of Simultaneous
             {D}iophantine Approximation Problems},
  journal = {SIAM Journal on Computing},
  volume  = {14},
  number  = {1},
  pages   = {196--209},
  year    = {1985},
  doi     = {10.1137/0214016}
}

@article{HochbaumMegiddoNaorTamir1993,
  author  = {Hochbaum, Dorit S. and Megiddo, Nimrod and
             Naor, Joseph and Tamir, Arie},
  title   = {Tight Bounds and 2-Approximation Algorithms for
             Integer Programs with Two Variables per Inequality},
  journal = {Mathematical Programming},
  volume  = {62},
  number  = {1--3},
  pages   = {69--83},
  year    = {1993},
  doi     = {10.1007/BF01585160}
}

@article{HochbaumNaor1994,
  author  = {Hochbaum, Dorit S. and Naor, Joseph},
  title   = {Simple and Fast Algorithms for Linear and Integer
             Programs with Two Variables per Inequality},
  journal = {SIAM Journal on Computing},
  volume  = {23},
  number  = {6},
  pages   = {1179--1192},
  year    = {1994},
  doi     = {10.1137/S0097539793251876}
}

@article{BarYehudaRawitz2001,
  author  = {Bar-Yehuda, Reuven and Rawitz, Dror},
  title   = {Efficient Algorithms for Integer Programs with Two
             Variables per Constraint},
  journal = {Algorithmica},
  volume  = {29},
  number  = {4},
  pages   = {595--609},
  year    = {2001},
  doi     = {10.1007/s004530010075}
}

@inproceedings{LahiriMusuvathi2005,
  author    = {Lahiri, Shuvendu K. and Musuvathi, Madanlal},
  title     = {An Efficient Decision Procedure for {UTVPI}
               Constraints},
  booktitle = {Frontiers of Combining Systems},
  series    = {Lecture Notes in Computer Science},
  volume    = {3717},
  pages     = {168--183},
  publisher = {Springer},
  year      = {2005},
  doi       = {10.1007/11559306_9}
}

@inproceedings{TamoriKimura2019,
  author    = {Tamori, Takuya and Kimura, Kei},
  title     = {A Fast Algorithm for Unbounded Monotone Integer
               Linear Systems with Two Variables per Inequality
               via Graph Decomposition},
  booktitle = {WALCOM: Algorithms and Computation},
  series    = {Lecture Notes in Computer Science},
  volume    = {11355},
  pages     = {209--218},
  publisher = {Springer},
  year      = {2019},
  doi       = {10.1007/978-3-030-10564-8_17}
}

@article{SubramaniWojciechowski2023,
  author  = {Subramani, K. and Wojciechowski, Piotr},
  title   = {Integer Feasibility and Refutations in {UTVPI}
             Constraints Using Bit-Scaling},
  journal = {Algorithmica},
  volume  = {85},
  number  = {2},
  pages   = {610--637},
  year    = {2023},
  doi     = {10.1007/s00453-022-01048-1}
}

@inproceedings{deVlas2024,
  author    = {de Vlas, Jorke M.},
  title     = {On the Complexity of Integer Programming with
               Fixed-Coefficient Scaling},
  booktitle = {30th International Conference on Principles and
               Practice of Constraint Programming (CP 2024)},
  series    = {Leibniz International Proceedings in Informatics
               (LIPIcs)},
  volume    = {307},
  pages     = {35:1--35:9},
  publisher = {Schloss Dagstuhl -- Leibniz-Zentrum f{\"u}r
               Informatik},
  year      = {2024},
  doi       = {10.4230/LIPIcs.CP.2024.35}
}

@article{Hochbaum2004,
  author  = {Hochbaum, Dorit S.},
  title   = {Monotonizing Linear Programs with up to Two
             Nonzeroes per Column},
  journal = {Operations Research Letters},
  volume  = {32},
  number  = {1},
  pages   = {49--58},
  year    = {2004},
  doi     = {10.1016/S0167-6377(03)00074-9}
}

@inproceedings{DadushKohNaturaOlverVegh2024,
  author    = {Dadush, Daniel and Koh, Zhuan Khye and
               Natura, Bento and Olver, Neil and V{\'e}gh, L{\'a}szl{\'o} A.},
  title     = {A Strongly Polynomial Algorithm for Linear Programs
               with at Most Two Nonzero Entries per Row or Column},
  booktitle = {Proceedings of the 56th Annual ACM Symposium on
               Theory of Computing},
  pages     = {1561--1572},
  publisher = {ACM},
  year      = {2024},
  doi       = {10.1145/3618260.3649764}
}

@article{PritchardChakrabarty2011,
  author  = {Pritchard, David and Chakrabarty, Deeparnab},
  title   = {Approximability of Sparse Integer Programs},
  journal = {Algorithmica},
  volume  = {61},
  number  = {1},
  pages   = {75--93},
  year    = {2011},
  doi     = {10.1007/s00453-010-9431-z}
}

@article{FioriniJoretWeltgeYuditsky2025,
  author  = {Fiorini, Samuel and Joret, Gwena{\"e}l and
             Weltge, Stefan and Yuditsky, Yelena},
  title   = {Integer Programs with Bounded Subdeterminants and
             Two Nonzeros per Row},
  journal = {Journal of the ACM},
  volume  = {72},
  number  = {1},
  articleno = {3},
  pages   = {3:1--3:50},
  year    = {2025},
  doi     = {10.1145/3695985}
}

@inproceedings{BrianskiLassotaPekarkovaPilipczukReuter2026,
  author    = {Bria{\'n}ski, Marcin and Lassota, Alexandra and
               Pek{\'a}rkov{\'a}, Krist{\'y}na and
               Pilipczuk, Micha{\l} and Reuter, Janina},
  title     = {On Integer Programs That Look Like Paths},
  booktitle = {Integer Programming and Combinatorial Optimization},
  series    = {Lecture Notes in Computer Science},
  volume    = {16588},
  pages     = {158--172},
  publisher = {Springer},
  year      = {2026},
  doi       = {10.1007/978-3-032-28691-8_11}
}

@inproceedings{Schaefer1978,
  author    = {Schaefer, Thomas J.},
  title     = {The Complexity of Satisfiability Problems},
  booktitle = {Proceedings of the 10th Annual ACM Symposium on
               Theory of Computing},
  pages     = {216--226},
  publisher = {ACM},
  year      = {1978},
  doi       = {10.1145/800133.804350}
}

@book{schrijver1986,
  title     = {Theory of Linear and Integer Programming},
  author    = {Schrijver, Alexander},
  year      = {1986},
  publisher = {John Wiley \& Sons},
  address   = {Chichester, West Sussex, England},
  series    = {Wiley-Interscience Series in Discrete Mathematics},
  isbn      = {0471908541}
}

@article{BlairJeroslow1982,
  author  = {Blair, Charles E. and Jeroslow, Robert G.},
  title   = {The Value Function of an Integer Program},
  journal = {Mathematical Programming},
  volume  = {23},
  number  = {1},
  pages   = {237--273},
  year    = {1982},
  doi     = {10.1007/BF01583794}
}

@article{BasuMartinRyanWang2019,
  author  = {Basu, Amitabh and Martin, Kipp and Ryan, Christopher Thomas
             and Wang, Guanyi},
  title   = {Mixed-Integer Linear Representability, Disjunctions,
             and {C}hv{\'a}tal Functions---Modeling Implications},
  journal = {Mathematics of Operations Research},
  volume  = {44},
  number  = {4},
  pages   = {1264--1285},
  year    = {2019},
  doi     = {10.1287/moor.2018.0967}
}

@article{BordeauxKatsirelosNarodytskaVardi2011,
  author  = {Bordeaux, Lucas and Katsirelos, George and
             Narodytska, Nina and Vardi, Moshe Y.},
  title   = {The Complexity of Integer Bound Propagation},
  journal = {Journal of Artificial Intelligence Research},
  volume  = {40},
  pages   = {657--676},
  year    = {2011},
  doi     = {10.1613/jair.3248}
}

@inproceedings{EibenGanianKnopOrdyniak2018,
  author    = {Eiben, Eduard and Ganian, Robert and
               Knop, Du{\v s}an and Ordyniak, Sebastian},
  title     = {Unary Integer Linear Programming with
               Structural Restrictions},
  booktitle = {Proceedings of the Twenty-Seventh International
               Joint Conference on Artificial Intelligence},
  pages     = {1284--1290},
  publisher = {International Joint Conferences on Artificial
               Intelligence Organization},
  year      = {2018},
  doi       = {10.24963/ijcai.2018/179}
}

@inproceedings{McFarland2026,
  author    = {McFarland, Caleb},
  title     = {Totally {$\Delta$}-Modular Tree Decompositions of
               Graphic Matrices for Integer Programming},
  booktitle = {52nd International Workshop on Graph-Theoretic
               Concepts in Computer Science (WG 2026)},
  series    = {Leibniz International Proceedings in Informatics
               (LIPIcs)},
  volume    = {376},
  pages     = {33:1--33:18},
  publisher = {Schloss Dagstuhl -- Leibniz-Zentrum f{\"u}r
               Informatik},
  year      = {2026},
  doi       = {10.4230/LIPIcs.WG.2026.33}
}

\end{document}